\documentclass[11pt,fleqn]{article}

\usepackage{amsmath}
\usepackage{amssymb}
\usepackage{times}
\usepackage{mathptm}
\usepackage{ifthen}
\usepackage{graphics}
\usepackage{color}
\usepackage{epsfig}
\usepackage[pdftex,colorlinks,bookmarks,pdfpagemode=UseOutlines,linkcolor=black,
pagecolor=black,urlcolor=black,citecolor=black,letterpaper,pageanchor=false]{hyperref}
\usepackage[margin=1in]{geometry}
\usepackage[T1]{fontenc}
\usepackage{titlesec}
\usepackage[authoryear,round]{natbib}
\usepackage{titlesec}
\usepackage{appendix}
\usepackage{wrapfig}

\long\def\comment#1{}

\newcommand{\G}{\Gamma}

\begin{document}
\title{Platoons  of connected vehicles can double throughput  in urban  roads\thanks{This research was supported by the California Department of Transportation UCTC Award 65A0529 TO041.  We thank Rene Sanchez of Sensys Networks for the data in Figure 2, and Alex A. Kurzhanskiy, Gabriel Gomes, Roberto Horowitz, Sam Coogan and others in the Berkeley Friday Arterial seminar for stimulating discussions. }
 }
\author{Jennie Lioris, Ramtin Pedarsani, Fatma Yildiz Tascikaraoglu  and Pravin Varaiya\footnote{
Jennie.lioris@cermics.enpc.fr, Ecole des Ponts ParisTech; ramtin@ece.ucsb.edu, Electrical and Computer Engineering, University of California, Santa Barbara; fayildiz@yildiz.edu.tr, Control and Automation Engineering,  Yildiz Technical University; varaiya@berkeley.edu, Electrical Engineering and Computer Scienc, University of California, Berkeley.}}
\maketitle

\begin{abstract} 
Intersections are the bottlenecks of the urban road system because an intersection's capacity  is only a fraction of  the maximum flows that the roads  connecting to the intersection can carry.  This capacity can be  increased if vehicles  cross the intersections in platoons rather than one by one as they do today.  Platoon formation is enabled by connected vehicle technology.   This paper  assesses the potential mobility benefits of platooning.  It argues that saturation flow rates, and hence intersection capacity, can be doubled or tripled by platooning.   The argument is supported by the analysis of three  queuing models and by the simulation  of a  road network with 16 intersections and  73 links.  The queuing analysis and the simulations reveal that a signalized network with fixed time control will support an increase in demand by a factor of (say) two or three if all  saturation flows are increased by the same factor, with no change in the control.   Furthermore, despite the  increased demand vehicles will experience the same delay and travel time.  The same scaling improvement is achieved when the fixed time control is replaced by the max pressure adaptive control.  Part of the capacity increase can alternatively be used to reduce queue lengths and the associated queuing delay by decreasing the cycle time.  Impediments to the control of connected vehicles to achieve platooning at intersections  appear to be small.
\end{abstract}

\section{Introduction} 
Connected vehicle technology (CVT) has aroused  interest in the academic community as well as in automobile and IT companies.  Researchers are exploring ways to achieve the ambitious information, mobility and safety goals announced by the USDOT ITS Joint Program Office (\cite{JPO_CV}):

Information: ``Vehicle-to-infrastructure (V2I) capabilities and anonymous information from passengers' wireless devices relayed through dedicated short-range communications (DSRC) and other wireless transmission media, has the potential to provide transportation agencies with dramatically improved real-time traffic, transit, and parking data, making it easier to manage transportation systems for maximum efficiency and minimum congestion.''

Safety: CVT can ``increase situational awareness and reduce or eliminate crashes through vehicle-to-vehicle (V2V) and V2I data transmission. V2V and V2I applications enable vehicles to inform drivers of roadway hazards and dangerous situations that they can't see through driver advisories, driver warnings, and vehicle and/or infrastructure controls.''

Mobility: CVT can ``identify, develop, and deploy applications that leverage the full potential of connected vehicles, travelers, and infrastructure to enhance current operational practices and transform future surface transportation systems management.''

The most important academic program on CVT has been  the Safety Pilot Model Deployment  (\cite{SPMD}). Its  objective was to assess ``a real-world deployment of V2V technology,'' involving 300 vehicles that could receive and send  safety messages (curve speed warning, emergency electronic brake light, and forward collision warning) and 2300 other vehicles, equipped with communications-only devices with no safety applications, serving as `target' vehicles.  The pilot demonstrated the capability to exchange messages using DSRC equipment from different suppliers.  The NHTSA report (\cite{V2VNhtsa}) describes the  enhanced safety that V2V communications potentially offers in many pre-crash scenarios.  \cite{V2VEE} provides a  very readable account of  the report's highlights.

Considerable research effort has been devoted to the formation and performance of wireless vehicular ad-hoc networks (VANETs) in a mobile environment.  There also is  research  on vehicle control using V2V communication to execute tasks such as  cooperative adaptive cruise control (CACC), merging, and safely crossing an intersection  (\cite{ploeg, kianfar,cacc-milanes,CACC-twente,cacc-hafner,Malikopoulos_CVT}). These references indicate the potential for automating some driving functions, with possible side benefits  of mobility and safety.  The papers of \cite{dresner_AI,huang_CV,barth_CV} present simulations of control algorithms that coordinate the movement of connected and autonomous vehicles through an intersection with no traffic signals as such, but with a sophisticated ``agent'' who schedules the use of the intersection.

Commercial effort in CVT seeks to stimulate and meet consumer demand to connect their cars to the internet and operate their cellphone, send a message or play a specific song hands-free, and remotely monitor if the car is exceeding a preset speed.    These commercial efforts are not primarily intended to promote mobility or safety.   Driver assistance technologies, including adaptive cruise control (ACC), lane-departure warning and self-driving cars, do seek to enhance safety, but they are vehicle-resident technologies that do not rely on V2V communication. Neither academic nor commercial research is directly concerned with enhancing mobility.  Lastly, major auto companies have publicized their  research activity in autonomous vehicles, but no scientific details have been published. 

This paper explores the use of CVT to directly increase the capacity of urban roads based on the simple idea that if CVT can enable organizing vehicles into platoons, the capacity of an intersection can be dramatically increased by  a factor of two to three.  \S \ref{sec-capacity} sketches a scenario  showing how platooning increases an intersection's capacity.  The capacity increase will make it possible to support an increase in demand in the same proportion, with no change in the signal control.  Alternatively, some of the capacity increase can be used to reduce queues and the associated delay at intersections.  An argument for this `free' capacity increase is based on an analysis of three queuing models in \S \ref{sec-predict} that predict the mobility benefits of platooning.  The simulation exercise in  \S \ref{sec-sim} lends plausibility to those predictions.   \S \ref{sec-conc} summarizes the main conclusions of this study.

\section{Intersection capacity} \label{sec-capacity}
The Highway Capacity Manual (\cite{HCM}) defines an intersection's  capacity as 
\begin{equation} \label{1}
\text{Capacity} = \sum_i s_i \frac{g_i}{T},
\end{equation}
in which $T$ is the cycle time and, for lane  group $i$,  $s_i$ is the saturation flow rate and $g_i/T$ is the effective green ratio.  HCM takes $s_i = N \times s_0 \times f$: $N$ is the number of lanes in the group,  $s_0$ is the base rate in vehicles per hour (vph), and $f$ is an  `adjustment factor' that accounts for the road geometry and nature of the traffic.  The base rate $s_0$ is obtained from a thought experiment: it is the maximum discharge rate from  an infinitely long queue of vehicles facing a permanently green signal.  HCM recommends  $s_0 = 1,900$ vph, although empirical estimates can be as low as 1,200 vph.  Note that $s_i \times (g_i /T)$ is  the rate  at which vehicles in queue in group $i$ can potentially be served by the intersection.  So we also call it the \textit{service rate} in a queuing model of this lane group.

\begin{figure}
\centering
\includegraphics[width=3.7in]{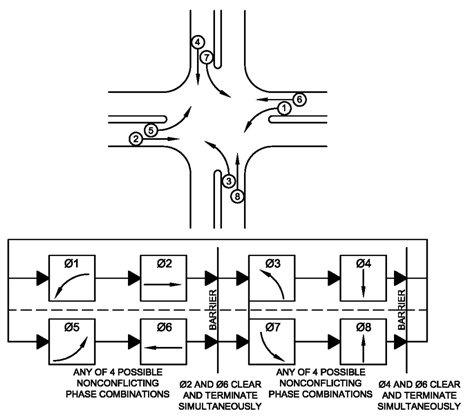}
\caption{An intersection with four approaches, each with a through and left-turn movement (top), and a ring-and-barrier diagram to organize the movements (bottom).  Source: \cite{Traffic_handbook}.}
\label{fig-fig1}
\end{figure}
Consider an intersection with four approaches, each  with one through lane and one left-turn lane as in Figure \ref{fig-fig1} (top).  There are thus eight movements in all.  Suppose each lane  supports a flow up to 1,900 vph for a total capacity of $1,900 \times 8 = 15,200$ vph.  However, from Figure \ref{fig-fig1} (bottom) only two movements can safely be allowed at the same time, so the effective green ratio for each movement is at most 0.25, and from equation \eqref{1} the capacity of the intersection is only 3,800 vph.   Thus the intersection is the principal bottleneck in urban roads: its  capacity is a  fraction (here 1/4) of the capacity of the roads connecting to it.  

\begin{figure}[h!]
\centering
\includegraphics[width=5.5in]{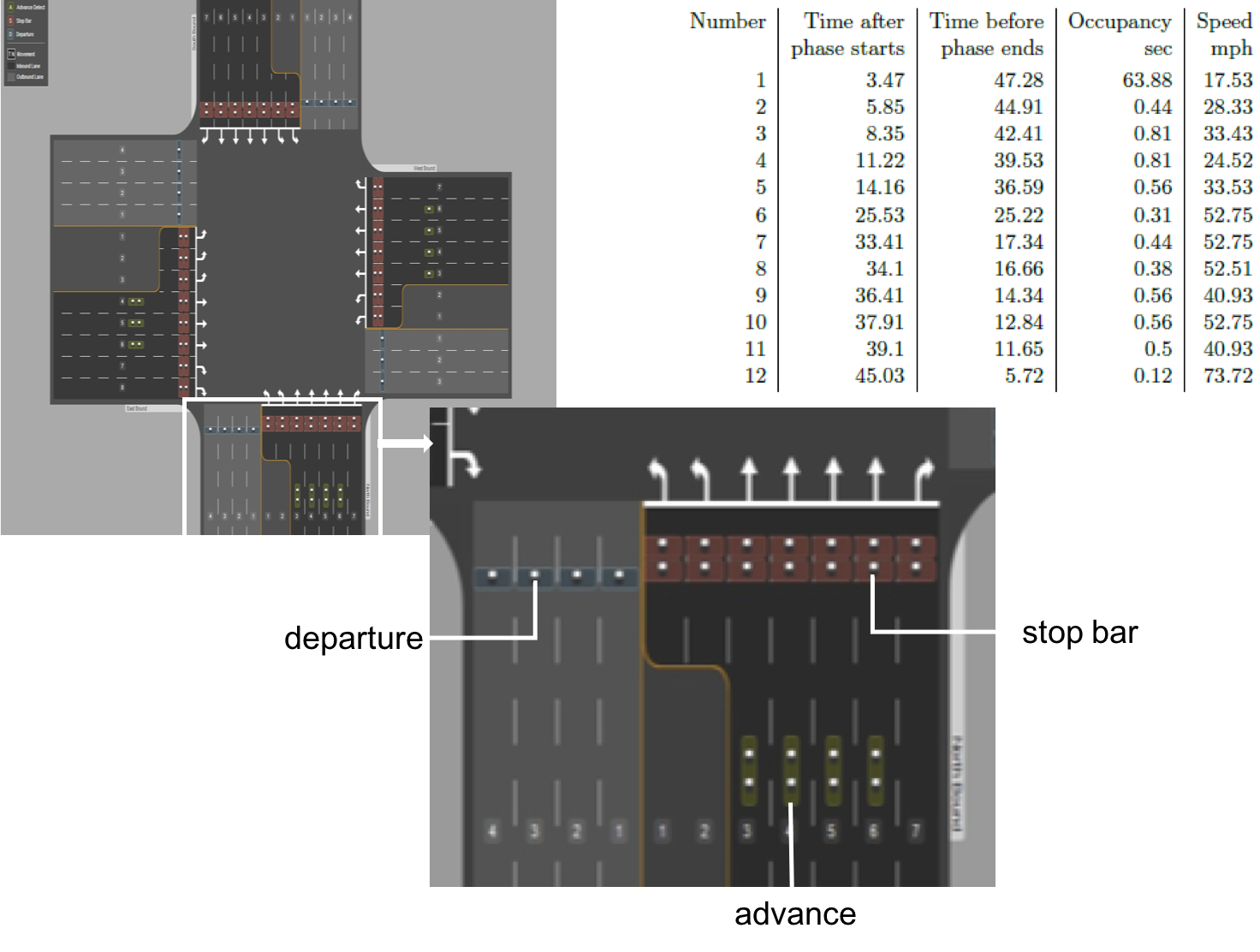}
\
\caption{Schematic of intersection detection system, and a trace of vehicles entering intersection from one through lane during one green phase.}
\label{fig:fig2}
\end{figure}
Figure \ref{fig:fig2} (left) is a schematic of the  system of vehicle detectors installed at an intersection in Santa Clarita, CA.  Each tiny white dot is a
magnetic sensor that reports the times at which a vehicle enters and leaves its detection zone.  When there is a pair of detectors, like at every stop bar and some advance  locations, the speed of the vehicle is also estimated and reported.  The detection system  also receives the signal phase timing from the conflict monitor card in the controller cabinet (not shown in the figure). The sensors in the departure lanes are used to estimate turn movements, as explained  in \cite{HRSensing}.
These measurements can be processed to obtain all intersection performance measures, including V/C (volume to capacity) ratios, fraction of arrivals in green, as well as red-light violations (\cite{bullock-performance,HR_TRC}).   From the times at which each vehicle enters the intersection we get the empirical saturation flow rate, that is the rate at which vehicles actually move through the intersection during the green phase, and the vehicle headway.

Figure \ref{fig:fig2} (top right) displays the trace of all 12 vehicles that enter a through lane in the intersection  during one cycle with a green phase duration of  50 sec for the through movement.  The second and third columns give the times (in sec) after the start of green and before the end of green when each vehicle enters the intersection, the fourth column gives the duration of time the vehicle `occupied' the detector zone, and the fifth column lists an estimate of its speed.  The detectors have a sampling frequency of 16 Hz, so the speeds are quantized (the distance between the two detectors is 12 feet), and speeds above 60 mph have a quantization error of about 15 mph, speeds below 30 mph have an error under 5 mph.  The average speed of the 12 vehicles is 42 mph.  The speed limit at this intersection is 50 mph.  The first vehicle entering the 
intersection has a delay or reaction time of 3.47 sec.  The first 5 vehicles enter within 14.16 sec, so the empirical saturation flow rate of this movement is
14.16/5 = 2.83 sec per veh or 3600/2.83 = 1272 veh/hour.  Vehicles 5, 6, $\cdots$ travel at much higher speed.

Suppose  these 12 vehicles were all to move as a platoon, that is, at the same speed of 45  mph (66 feet/sec) with a uniformly small time headway of (say) 0.75 s  (the space headway would be $0.75 \times 66 = 49.5$ ft), giving a saturation flow rate of $3600/0.75 = 4800$ vph, which is 3.8 times the observed rate of 1272 vph and 2.5 times HCM's theoretical rate of 1900 vph.  
The headways among vehicles 1, $\cdots, 5$ are small, suggesting that they were queued at the stop bar.
The larger headway between vehicles 6 and 7 and 7 and 8 suggests they were moving when the phase started.
The small headway among vehicles $7, \cdots, 11$ suggest they were moving as a platoon. An estimate for an intersection with speed of 30 mph (44 feet/sec) and a space headway of 40 feet is a (platoon) saturation flow rate of (40/44) $\times$ 3600 = 3960 vph, which is twice HCM's  1900 vph and up to  three times the  rates observed in today's intersections. The following fact summarizes the relationship between platooning and the increase in saturation flow rates.

\textbf{Fact} 
Platooning  decreases the headway of vehicles crossing an intersection by a gain factor $\G$ and increases the saturation flow rates of an intersection by the same factor $\G$.

What does it take to organize a platoon? If each of the 12 vehicles queued at (or approaching) the intersection in Figure \ref{fig:fig2} could measure, say by radar, the relative distance and speed from the vehicle in front of it, its longitudinal motion could then be 
controlled by an adaptive cruise control (ACC) algorithm that would maintain a tight headway.  If furthermore, these vehicles could communicate with each other and with the signal controller (for phase information), a cooperative adaptive cruise control (CACC) algorithm would maintain an even shorter headway and thus achieve a greater saturation flow rate.

Platooning is technically feasible.  In fact, an autonomous 8-vehicle platoons with 16 ft gap traveling at 60 mph was demonstrated 20 years ago in 1997 by the National Automated Highway System Consortium  (\cite{NAHSRC, shladover-keynote}).  Since then vehicle automation has been greatly facilitated by advances in actuation (electronic braking, throttle, and steering) and sensing (radar and video), while platoon stability and control design are much better understood.  
Indeed, \cite{ploeg} report an experiment of a 6-vehicle CACC platoon, with a 0.5s headway.  They used ACC equipped vehicles  augmented with V2V communications using a  802.11a WiFi radio in ad-hoc mode.  \cite{cacc-milanes} describe experiments with a 4-vehicle platoon, capable of cut-in, cut-out and other maneuvers, using CACC technology.  The vehicles'  factory-equipped ACC capability was enhanced by a DSRC radio that permitted V2V communication  to enable CACC. The vehicles in the platoon  had a time gap of 0.6 sec (time headway of 0.8 s) traveling at 55 mph.  ACC is today common in many high-end cars.

Of course challenges in control and communication will be encountered and must be overcome, if the  ideas proposed here and in the cited references are to be realized.  However the challenges are much less difficult for our requirement that (C)ACC mode must be enabled only in vehicles that are stopped at intersections so that they can  cross the intersections as platoons.  The paper explores what throughput and delay advantages such organization can bring.

It is important to distinguish our proposal to use (C)ACC to increase an intersection's capacity from proposals to use CACC to increase  road capacity by decreasing headway. Increasing the capacity of urban roads will \textit{not} increase the throughput of the urban network which is limited by intersection capacity.

In this paper we consider the limiting case in which all vehicles are  connected.  But for several years into the future there will be a mixture of regular and connected vehicles, so the increase in the saturation flow rate will depend on the market penetration rate.  However, our mathematical models and  simulations  predict   performance when saturation flows are (on average) increased by \textit{any} factor $\G$. Of course when there are fewer connected vehicles, this factor $\G$ will be smaller. In the  12-vehicle  queue of the example above, suppose vehicles 2 and 5 are manually driven. Suppose that a platoon must start with a connected vehicle. Then  the \textit{average}  headway will be $\bar{H}= (4 H_{high} + 7 H_{low})/11$, where $H_{high}$ is the headway of regular vehicles and $H_{low}$ is the headway of connected vehicles. For $H_{high} = 2$ sec and $H_{low} = 0.75$ sec, one finds $\bar{H} = 1.2$ which results in $\G \simeq 1.67$, instead of $\G \simeq 2.67$ if all vehicles are connected. Of course, the exact value of $\G$ depends on the setup and the technology that is used, and evaluating their effect is  an interesting  research direction.

\section{Three predictions}\label{sec-predict}
Suppose magically that the saturation flow rates at all signalized intersections in an urban road network are increased by the same factor of 2 to 3 using platoons.  It  seems intuitively obvious that this change  will increase throughput and reduce travel time.  It is less clear how much the improvement would be, since the  links at an intersection may get filled up by vehicles arriving more rapidly than at the normal rate, and thereby block additional arrivals.

An exact analysis of the stochastic queuing network that models the road network is beyond our reach.   So we use three   models to predict the average queue size and delay of vehicles, and   throughput, beginning with the simplest  model for an isolated intersection, namely the M/M/1 queue (\cite{soh, anokye, pavone}).  We will see that the predictions from the M/M/1 queue are unduly optimistic because its service process does not capture the stop-go nature of signal lights. 

For the second model we propose a memoryless queue with an on/off (green/red) service process that  better models signal actuation.  We use the Recursive Renewal Reward (RRR)   technique to exactly analyze the Markov chain corresponding to this queue, and characterize the  behavior of the system as both arrival  and saturation rates increase in the same proportion. The predictions of this second model are more reasonable.

Our third model is a fluid queuing network proposed in \cite{FTControl}, and we study its behavior as saturation flows and arrival rates are increased by the same factor.  We show consistency between the  predictions of the second model and this fluid  network model.  We develop an example to indicate the challenge in dealing with  finite capacity queues. 

\subsection{Model 1: M/M/1 queue}
Consider an M/M/1 queue with arrival rate $\lambda$ and service rate $\mu$ as the model of the queue at a single approach of an isolated intersection.  Vehicles arrive at this approach at rate $\lambda$ vph and are served at rate $\mu$ vph, the saturation flow rate multiplied by the effective green ratio of the approach (see \eqref{1}).

As is  known  the stability region of the system is $\lambda < \mu$, and the average delay or time spent in the queue  by a vehicle is $\frac{1}{\mu}\times\frac{\rho}{1 - \rho}$, where $\rho = \lambda / \mu$ is the load of the system ($\rho$ is also the `volume-capacity or VC ratio'). Now if $\mu$ is increased to $\G \mu$ for $\G > 1$, the stability region is increased $\G$-fold to $\lambda < \G \mu$. Further, if the demand is also increased to $\G\lambda$, the average delay of the system is \textit{decreased} by a factor $\G$ since the load $\rho$ is unchanged. This  suggests that   increasing the saturation flow rates by a factor $\G$ leads to an increase in  throughput and a decrease in delay by the same factor $\G$. This double benefit is too good to be true. 

One immediate objection to the `double benefit' is that the queues at an intersection have finite capacity so it is more appropriate to consider the M/M/1/$K$ queue, where $K$ is the finite capacity of the queue. Again for $\rho = \frac{\lambda}{\mu}$ the stationary distribution of the length of the M/M/1/$K$ queue  is 
\begin{align} \label{2}
\pi_k = \frac{\rho^k}{\sum_{i=0}^K \rho^i} = \rho^k \frac{1 - \rho}{1 - \rho^{K+1}},\; k = 0,1,\ldots, K,
\end{align}
in which $\pi_k$ is the probability that there are $k$ vehicles in queue.  New arrivals are blocked when the queue is full, which occurs with probability 
$\pi_K$; thus, the effective throughput of the system is $\lambda_e = \lambda (1 - \pi_K)$. The expected queue length is
\begin{align}\label{3}
\bar{N} = \sum_{k=0}^K k \pi_k = \frac{\rho(1 - (K+1)\rho^K + K \rho^{K+1})}{(1-\rho)(1 - \rho^{K+1})},
\end{align}
and the average delay by Little's law is $\bar{D} = \frac{\bar{N}}{\lambda_e}$. 

As one can see from the formulas, the  blocking probability $\pi_K$ and hence the throughput  depends only on $\rho$. Thus if $\mu$ and $\lambda$ are increased by the same factor $\G$,  the  blocking probability  $\pi_K$ is unchanged.  Hence the   throughput will be  increased by  $\G$, since $\lambda_e = \lambda (1 - \pi_K)$, while the delay will be decreased by 
$\G$.  (The blocking probability is small for reasonable parameters: with $\rho$ or VC ratio equal to 0.8 and $K=10$, the blocking probability is 0.02.)  Thus the double benefit of the M/M/1 model is assured even with finite capacity links.

The M/M/1/$K$ Model  confirms one intuition: if the saturation flow of all links along an arterial is increased by factor $\G$, it is the shortest link (smallest $K$) that will limit the increase in throughput.  But this intuition may  be too simplistic in ignoring the possibility of designing offsets so that platoons travel in a green wave unhindered by the small storage capacity of short links rather than in the Poisson stream of the M/M/1/$K$ model.  However, the main reason the M/M/1/$K$ model is unduly optimistic is that it replaces the on-off (green-red) nature of signal actuation by a constant saturation rate.  The next model does not have this defect.

\subsection{Model 2: Single queue with on/off memoryless service}\label{sec:single}
\begin{figure}[h!]
\centering
\includegraphics[scale = 0.5]{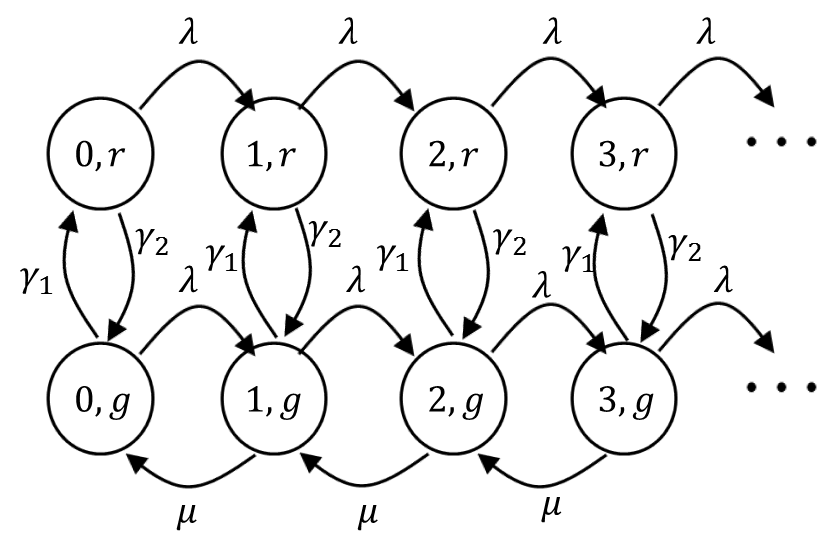}
\caption{Markov chain of a memoryless queue with on/off service.}
\label{fig:mc}
\end{figure}
Consider a single queue with on/off service: it receives full service at the saturation flow rate, $\mu$, when the light is green and no service when the light is red.  The queue behavior is modeled as the continuous-time Markov chain of Figure \ref{fig:mc} with states denoted $(x,l)$ in which $x$ is the queue length and $l$ is the state of service, with $l =g$ when there is  service (light is green), and  $l=r$ when there is no service (light is red).\footnote{The context will permit the distinction between the $g$ used for green duration in \eqref{1}, \eqref{14b}, the $g$ used for green light in Figure \ref{fig:mc}, and the $g$ used
for speedup factor in \eqref{11a}, \eqref{15a}.}

The state transition rates are given by the labels. 
This means, for example, that the transition probabilities from state $(x,l) = (0,r)$ are  given by:
\begin{eqnarray*}
P\{(x,l)(t+dt) = (1,r)~|~(x,l)(t) = (0,r)\} &=& \lambda dt +o(t),\\
P\{(x,l)(t+dt) = (0,g)~|~(x,l)(t) = (0,r)\} &=& \gamma_2 dt +o(t) ,\\
P\{(x,l)(t+dt) = (0,r)~|~(x,l)(t) = (0,r)\} &=& 1 - \lambda dt - \gamma_2 dt +o(t).
\end{eqnarray*}
We show that like in the M/M/1 model there is a throughput benefit but unlike in that model there is no delay benefit when the demand and service rates ($\lambda$ and $\mu$) both are  increased by  $\G$.

We compute the average delay corresponding to the Markov chain in Figure \ref{fig:mc}  using the \emph{Recursive Renewal Reward} (RRR) technique  described by \cite{gandhi}. RRR is based on renewal reward theory and busy period analysis of a 2-dimensional Markov chain like the one in Figure \ref{fig:mc} that is finite in one of the dimensions. This Markov chain can also be analyzed by matrix-analytic methods (\cite{latouche}), but  the RRR technique is better as it does not require computing the stationary distribution of the  chain  to determine the average delay (or other performance metrics). 

First note that the chain is in a green state for the fraction $\frac{\gamma_2}{\gamma_1 + \gamma_2} $ of the time on average, so in accordance with \eqref{1} the effective service rate is $\frac{\gamma_2}{\gamma_1 + \gamma_2}  \mu$ and the chain is stable if $\lambda < \frac{\gamma_2}{\gamma_1 + \gamma_2} \mu$.  We define a renewal cycle as  the process that starts from state $(0,g)$ and returns back to this state for the first time.  Imagine  that the process collects a reward in each state $(x, l)$ equal to the number $x$ of vehicles in the queue $(l = g, r)$. 
Let $R$ be the expected  reward in a cycle, and let $T$ be the expected duration of a cycle. Then the average number of vehicles in the queue is $R/T$, so the average delay is $\frac{R}{\lambda T}$. We use the birth-death property of the Markov chain to find linear recursive equations to calculate $R$ and $T$ as follows. Let $T_{x,g}^L$ be the expected time to visit state $(x-1,g)$ from state $(x,g)$ (expected time to go one level left), and let $T_{x,r}^L$ be the expected time to visit state $(x,g)$ from state $(x,r)$. Let $R_{x,g}^L$ be the expected  reward collected starting from state $(x,g)$  to state $(x-1,g)$, and let $R_{x,r}^L$ be the expected reward accumulated starting from state $(x,g)$ to state $(x,r)$. 
We then have the following linear relations:
\begin{align}
T &= \frac{1}{\lambda + \gamma_1} + \frac{\lambda}{\lambda + \gamma_1} T^L_{1,g} + \frac{\gamma_1}{\lambda + \gamma_1} T_{0,r}^L ,\label{4}\\
T_{1,g}^L &= \frac{1}{\lambda + \mu + \gamma_1} + \frac{\lambda}{\lambda + \mu + \gamma_1}(T_{1,g}^L + T_{2,g}^L) + \frac{\gamma_1}{\lambda + \mu + \gamma_1}(T^{L}_{1,r} + T_{1,g}^L), \label{5}\\
T_{0,r}^L &= \frac{1}{\lambda + \gamma_2} + \frac{\lambda}{\lambda + \gamma_2}(T^L_{1,r} + T^L_{1,g}). \label{6}
\end{align}
By the repetitive structure of the Markov chain, we have $T^L_{1,r} = T^L_{0,r}$ and $T_{2,g}^L = T_{1,g}^L$. Thus, one can solve these three linear equations to find the three `unknowns' $T$, $T^L_{1,g}$ and $T^L_{0,r}$. One can similarly write linear relations for the reward: 
\begin{align}
R &= \frac{0}{\lambda + \gamma_1} + \frac{\lambda}{\lambda + \gamma_1} R^L_{1,g} + \frac{\gamma_1}{\lambda + \gamma_1} R_{0,r}^L ,\label{7}\\
R_{1,g}^L &= \frac{1}{\lambda + \mu + \gamma_1} + \frac{\lambda}{\lambda + \mu + \gamma_1}(R_{1,g}^L + R_{2,g}^L) + \frac{\gamma_1}{\lambda + \mu + \gamma_1}(R^{L}_{1,r} + R_{1,g}^L) ,\label{8}\\
R_{0,r}^L &= \frac{0}{\lambda + \gamma_2} + \frac{\lambda}{\lambda + \gamma_2}(R^L_{1,r} + R^L_{1,g}). \label{9}
\end{align}
By Theorem 2 of \cite{gandhi}, since the reward at each state is the number of vehicles in the queue, $R_{2,g}^L = R_{1,g}^L + T^{L}_{1,g}$ and $R^L_{1,r} = R^L_{0,r} + T^L_{0,r}$. Solving these equations one finds $R$, $R^L_{1,g}$, and $R^L_{0,r}$.  The average queue length is
\begin{align}
\bar{N} = \frac{R}{T} = \frac{\lambda \gamma_1^2 + 2 \lambda \gamma_1 \gamma_2 + \lambda \gamma_1 \mu + \lambda \gamma_2^2}{(\gamma_1 + \gamma_2) (\gamma_2\mu - \lambda (\gamma_1 + \gamma_2))}. \label{10}
\end{align}
As a  check,  observe that when $\gamma_1 = 0$, the expression reduces to $\lambda/(\mu - \lambda)$, the same as for the M/M/1 queue, as  expected. Further, when $\lambda$ approaches the capacity limit $\frac{\gamma_2}{\gamma_1 + \gamma_2}\mu$, the expected number of vehicles in the queue tends to infinity. By Little's law, the average delay is 
\begin{align}
\bar{D} = \frac{\bar{N}}{\lambda} = \frac{\gamma_1 + \gamma_2 + \frac{\gamma_1}{\gamma_2 + \gamma_1} \mu}{\gamma_2 \mu - \lambda (\gamma_1 + \gamma_2)}.\label{11}
\end{align}    
We study how the delay and queue length change as both demand and service rates and cycle times are scaled.  To focus ideas, take the base case to be $\mu^0 = 2000$ vph, $\lambda^0 = 900$ vph, $\gamma_1^0 = \gamma_2^0 = \gamma^0 = 30$ switches per hour, corresponding to 30 cycles per hour or a cycle time of 120 sec and an effective green ratio of 0.5.  We scale these parameters as $\mu = \G \mu^0, \lambda = \G \lambda^0$, $\gamma_1 = \gamma_2 = g \gamma^0$.   Substituting these values into \eqref{11} gives
\begin{equation}\label{11a}
\bar{D} = \frac{1}{g} \frac{60g + 1000 \G}{30 \times 200 \G}, \; \bar{N} = \G  \lambda^0 \times \bar{D}.
\end{equation}
For  $g=1$ (no change in cycle time) we see that the average delay $\bar{D}$ does not materially change as $\G$ is increased to 2 or 3, whereas the average queue length $\bar{N}$ increases linearly with $\G$.  Thus this on-off model suggests two predictions as the saturation flow rate increases by a factor $\G$: (1) the network can support a demand that increases
by the same factor, (2) the queue length ($\bar{N}$) grows by the same factor, but the queuing delay ($\bar{D}$) is unchanged.
Another prediction can be extracted by considering the effect of an increase in $g$, which decreases the cycle time by the factor $g$, while keeping the effective green ratios the same.  We see from  expression \eqref{11a} that the average delay $\bar{D}$ and the queue size $\bar{N}$ \textit{decrease} almost linearly with $g$.  This gives a third prediction: (3) the  queue length can be or reduced by reducing the cycle time.  

For an isolated intersection a simple intuition supports these three predictions.  We assume that the chain is stable so that the queue builds up during
red and discharges during green.  The build up of the queue will be proportional to the product of the arrival rate and the red duration, so if the arrival rate is increased by factor $\G$, the queue size will also increase by the same factor.  On the other hand, if the discharge or service rate is increased by $\G$, the queue will be drained $\G$ times more quickly, so the delay will not change.  Now the red duration is proportional to the cycle time so if the cycle time is reduced by a factor $g$, the queue length will be reduced by the same proportion $g$.  

This intuition fails when we consider a network of intersections, since
a more rapid discharge of one queue at one intersection will lead to an equally rapid increase in arrivals at a queue in the next intersection. Furthermore, in a network model we need to explicitly account for the travel time over links and the effect of  phase offsets on queue formation.
It is remarkable, therefore, that these same predictions, suitably reformulated, hold for a network as we study next.

\subsection{Model 3: Fluid model of queuing network} \label{sec-fluid}

\begin{figure}[h!] 
\centering
\includegraphics[width=4in]{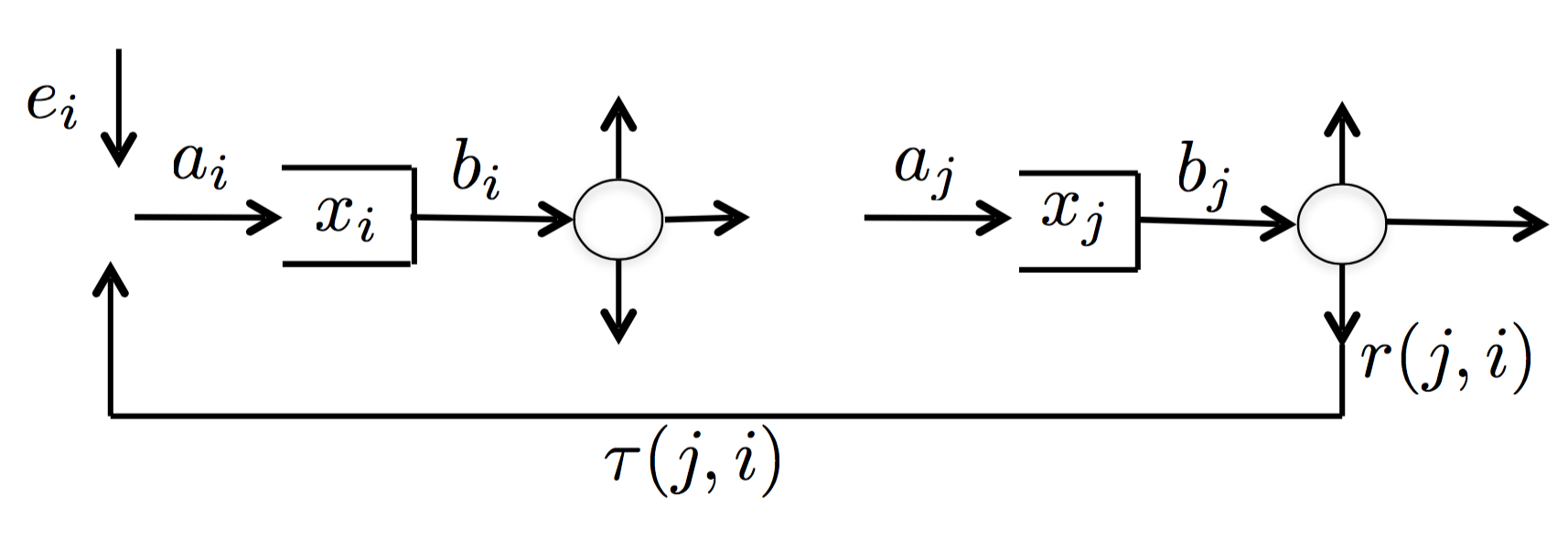}
\caption{Network of fluid model.  Source: \cite{FTControl}.}
\label{fig:fig4}
\end{figure}
We now analyze a  fluid model of a network of signalized intersections, all with fixed-time  control with the same cycle time $T$. The  model was studied by \cite{FTControl} and illustrated in Figure \ref{fig:fig4}. There are $J$ queues in the network, each corresponding to one turn movement.  A specified fraction $r(j,i)$ of the vehicles departing from queue $j$  is routed  to  queue $i$. 
Queue $i$ also has  exogenous  arrivals with rate $e_i(t)$. 
The  exogenously specified service rate of queue $i$, $c_i(t)$, is  periodic with period $T$: $c_i(t)$ is the saturation flow rate when the light is green and $c_i(t) = 0$, when it is red.   
$x_i(t)$ is the queue length of $i$  and $b_i(t)$ is its departure rate at time $t$.  $a_i(t)$ denotes the total arrival rate into queue $i$ at time $t$. The travel time from queue $j$ to queue $i$ on link $(j,i)$ is $\tau(j,i)$ as in the figure. The network dynamics are  as follows. 
\begin{align} 
\dot{x}_i(t)& = a_i(t) - b_i(t) , \label{12} \\
a_i(t) &= e_i(t) + \sum_{j=1}^J b_j(t-\tau(j,i)) r(j,i) , \label{13}\\
b_i(t) &= 
\left \{
\begin{array}{ll}
c_i (t) , & \mbox{ if } x_i(t) > 0, \\
\in [0, c_i (t)], & \mbox{ if } x_i(t) = 0,\\
0, & \mbox { if } x_i(t) < 0.
\end{array}
\right . \label{14}
\end{align}    
\textbf{Remark} 
\begin{figure}[h!] 
\centering
\includegraphics[width=4in]{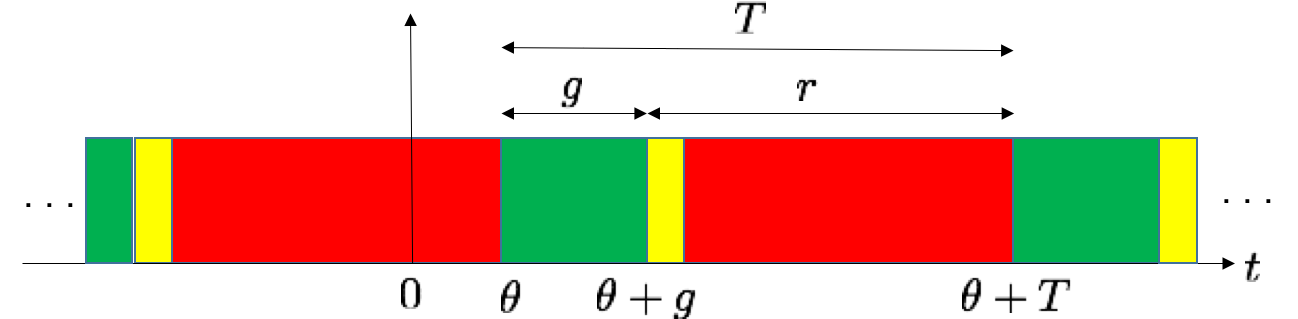}
\caption{The timing  for a movement whose start of green has offset $\theta$, green duration $g$, and cycle time $T$. }
\label{fig-figs}
\end{figure}
We briefly describe how to construct the service rate $c_i(t)$ from a timing diagram and offset values.  Figure \ref{fig-figs} is a timing diagram of a particular movement with  start of green offset by $\theta$ (relative to $t=0$ of a master clock), green duration $g$, red duration (including yellow) $r$, and cycle time $T = g + r$.  The
timing diagram is periodic.  Suppose the saturation rate for this movement is $C$.  Then the service rate is, for any integer $n$, 
\begin{equation}\label{14b}
c(t) = \left \{
\begin{array}{ll}
C, & t \in [nT + \theta, nT + \theta +g] \mbox{  -- the green intervals}\\
0, & t \in [nT + \theta +g, \theta + (n+1)T] \mbox{  -- the red intervals}
\end{array}
\right .
\end{equation}
If the signal were actuated or adaptive, the right-hand side of \eqref{14b} would be replaced by an expression for the actuation interval, which would no longer be periodic.  An example of an adaptive signal control is the `max pressure' controller introduced in section \ref{sec-sim} and described in 
\cite{MP,LiorisTRR,wodes}.

Returning to \eqref{12}-\eqref{14}, note that this is a nonlinear differential inclusion  with delay: nonlinear because $b_i(t)$ is nonlinear in $x_i(t)$, inclusion (rather than equation) because the right-hand side 
of  \eqref{14} is set-valued, and the $\tau(j,i)$ introduce delay.  The inclusion is not Lipschitz, hence neither  existence nor uniqueness of a solution is assured.  However, a surprising result of 
\cite{FTControl} guarantees that this differential inclusion has a unique solution for $x_i(0) \ge 0$.
We regard as exogenous parameters of the system  the rates $\{e_i(t), c_i(t)\}$ and the initial state $x(0)$.  

The system \eqref{12}-\eqref{14} is positively homogeneous of degree one.   That is, if $\{x(t), a(t), b(t)\}$ is the solution for a specified $\{x(0), e(t), c(t)\}$ and $\G > 0$, then $\{\G x(t), \G a(t), \G b(t)\}$ is the solution for $\{\G x(0), \G e(t), \G c(t)\}$.  Hence if the exogenous arrivals $e$ and the saturation flow rate are both increased by a factor $\G$, the  queue $x(t)$ and the departures $b(t)$  both increase by the same factor.  Thus the model predicts the first benefit of higher service rate (by platooning), namely, throughput increase in the same proportion.

What about queue length and delay?  The average queue length at $i$ over (say) one cycle beginning at time $t_0$ in the base case ($\G =1$) is
\[\bar{x}_i (t_0)= \frac{1}{T} \int_{t_0}^{t_0 + T} x_i(t) dt,\]
and if the service  and arrival rates (and the initial queue $x(0)$) are increased by a factor $\G$, the average queue length over the same cycle
also increases by the same factor:
\begin{equation}
\frac{1}{T} \int_{t_0}^{t_0 + T} \G x_i(t) dt = \G \bar{x}_i (t_0).\label{14a}
\end{equation}
Since the queue length, arrival and departure processes all increase by the same factor, it is immediate that the average delay per vehicle in each queue is \textit{unchanged}.
Since the travel time is the sum of the queuing delays and the free flow travel time along a route, and both of these are unchanged by the scale factor $\G$, it follows that the travel time is unchanged as well.  This result does not even require  the system to  be stable, i.e. queues to be bounded.  

In summary, the on-off queue (Model 2) and the fluid network (Model 3) yield a consistent prediction: If the saturation flow for every movement is increased by a factor $\G$, the network can support a throughput that increases by $\G$, while keeping the queuing delay and travel time  unchanged.

\subsection{Queues with finite capacity: An example}
In all three models above, the queue capacity is taken to be infinite.  Since the queues in the second and third models grow in proportion to the scale factor $\G$ of the saturation rates, the link may become saturated, and the throughput gain may be less than $\G$.  We illustrate this in a simple example using the fluid model.  Rigorous analysis of the queuing network  with finite capacity queues is difficult. 

Consider a single queue with  capacity of 20 vehicles, so arrivals are blocked once the queue length reaches 20.  Consider a constant arrival rate, $e(t) = a(t) = 10$, and the following service rate in one period $T = 2$: 
\begin{equation*}
c(t) = \left \{ \begin{array}{ll}
0, & \quad t \in [0,1) \\
30, & \quad t \in [1,2)
\end{array}
\right . .
\end{equation*}
Then $\bar{c} = 15 > \bar{a} = 10$, so the queue is stable. If $x(0) = 0$, the trajectory  is periodic with period $T$ and easily calculated to be
\begin{equation*}
x(t) = \left \{ \begin{array}{ll}
10t, & \quad t \in [0,1) \\
\max(30 - 20t,0), & \quad t \in [1,2)
\end{array}
\right . .
\end{equation*}
The maximum queue-length occurs at the end of red, at $t = 1 + nT, ~ n=0, 1, \cdots$, and equals $x(1) = 10$. 

Now suppose that the saturation flow rate and the arrival rate are increased by a factor of 3. Then $\bar{c} = 3 \times 15 = 45$ will also be the maximum throughput of the system in the case of infinite capacity, and the queue will increase by a factor of 3, so the maximum queue length will be $30$. This is larger than the 20 vehicle limit for the finite capacity queue. If the arrival rate is $a(t) = 3 \times 10 = 30$, the queue will be blocked during the interval $t \in [2/3,1]$. 
The maximum throughput of the finite-capacity queue is easily calculated to be $\frac{20 + 30}{2} = 25 = 5/9 \times 45$, which is only 55 percent of the  throughput of the  infinite-capacity queue.

The storage capacity of the link  rather than the intersection  has now become the bottleneck, so the throughput gain is reduced  to a  sub-linear function of the gain $\G$ in the saturation flow rates.  Note, however, that in a network context, whether a link becomes a bottleneck is difficult to determine  since it depends on the signal control and the offsets.

\textbf{Remark} Suppose the cycle length is decreased to $T = 2/3$ and the green ratio is  unchanged   while the saturation flow rate is increased by 3.  Then the throughput will still increase by $3$, as no blocking of the queue occurs, and again, it is the intersection that is the bottleneck.   This is an instance of prediction (3): queues are reduced if the cycle length is reduced, while keeping the green ratios (the $g_i /T$) in \eqref{1} the same.  The prediction will be validated in the simulations of the next section.
There is another moral: if the cycle length is reduced more and more, the M/M/1 model becomes a better fit and the second benefit of lower  delay begins to appear.

\subsection{Reducing cycle time, queue length and delay}\label{sec-reduce}
We formalize the intuition in the previous example.  In the model
\eqref{12}-\eqref{14} change the service rate to $c(gt)$, where $g > 1$.  This means that the cycle time is reduced by $g$, but the green ratios are unchanged.  Let $z(t) = g^{-1}x(gt)$.  Then
\[\dot{z}_i(t) = \dot{x}_i (gt) =a_i(gt) - b_i(gt),\]
\[
b_i(gt) = 
\left \{
\begin{array}{ll}
c_i (gt) , & \mbox{ if } x_i(gt) > 0 \leftrightarrow z_i(t) > 0, \\
\in [0, c_i(gt) (t)], & \mbox{ if } x_i(gt) = 0 \leftrightarrow z_i(t) = 0,\\
0, & \mbox { if } x_i(gt) < 0 \leftrightarrow z_i(t) < 0.
\end{array}
\right . 
\]
Furthermore from \eqref{13}
\[a_i(gt) = e_i(gt) + \sum_{j=1}^J b_i(gt - \tau(j,i))r(j,i).\]
Hence 
\begin{equation}\label{15a}
z(t) = g^{-1}x(gt)
\end{equation}
 is the queue for exogenous arrivals $e(gt)$, and service $c(gt)$.  Thus by speeding up the service rate by factor of $g$, the queue length is reduced by the same factor $g$.  This is the third prediction.

The third prediction is of a different character, since it has nothing directly to do with increasing saturation flow rate.  However, in practice reducing cycle time is not possible because it will lead to a reduction in intersection capacity (the constant lost time due to the yellow and all-red clearance intervals reduces the effective green ratio in \eqref{1}), and the base case demand may not be accommodated.  But if platooning gives a gain of $\G$ in the saturation flow rate, some of the gain may be used to reduce cycle time, queue length and delay.

\section{Case study}\label{sec-sim}

We now present a simulation study of a road network near Los Angeles, using a mesoscopic simulator called PointQ.  The simulator and the network are described in \cite{pointqa,wodes}, together with a base case  with exogenous demands at the input links, modeled as stationary Poisson streams, and intersections  regulated by fixed time (FT) controls and offsets.  
PointQ is a discrete event simulation; it accurately models vehicle arrivals, departures and signal actuation.  It models queues as `vertical' or `point' queues that discharge at the saturation flow rate when the signal is green,
and  are fed by exogenous  arrivals or by vehicles that are routed from other queues.  When a vehicle is discharged from one queue it travels to a randomly assigned destination queue according to the probability distribution specified by the routing matrix $ \{r(i,j)\}$ (see \eqref{13}).  The vehicle takes a pre-specified fixed time to travel along the link determined by the assigned destination and then joins the destination queue.  Every event in the simulation is recorded and uploaded into a database from which the reported performance measures are calculated.  

\begin{figure}[h!]
\centering
\includegraphics[width=6in]{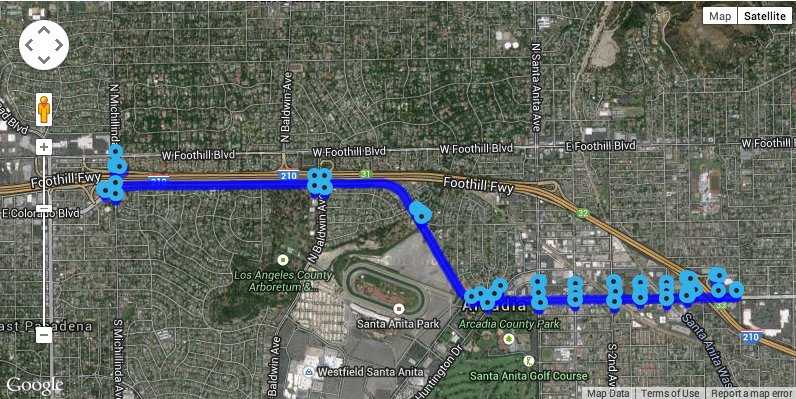}

\vspace*{0.5in}

\includegraphics[width=6in]{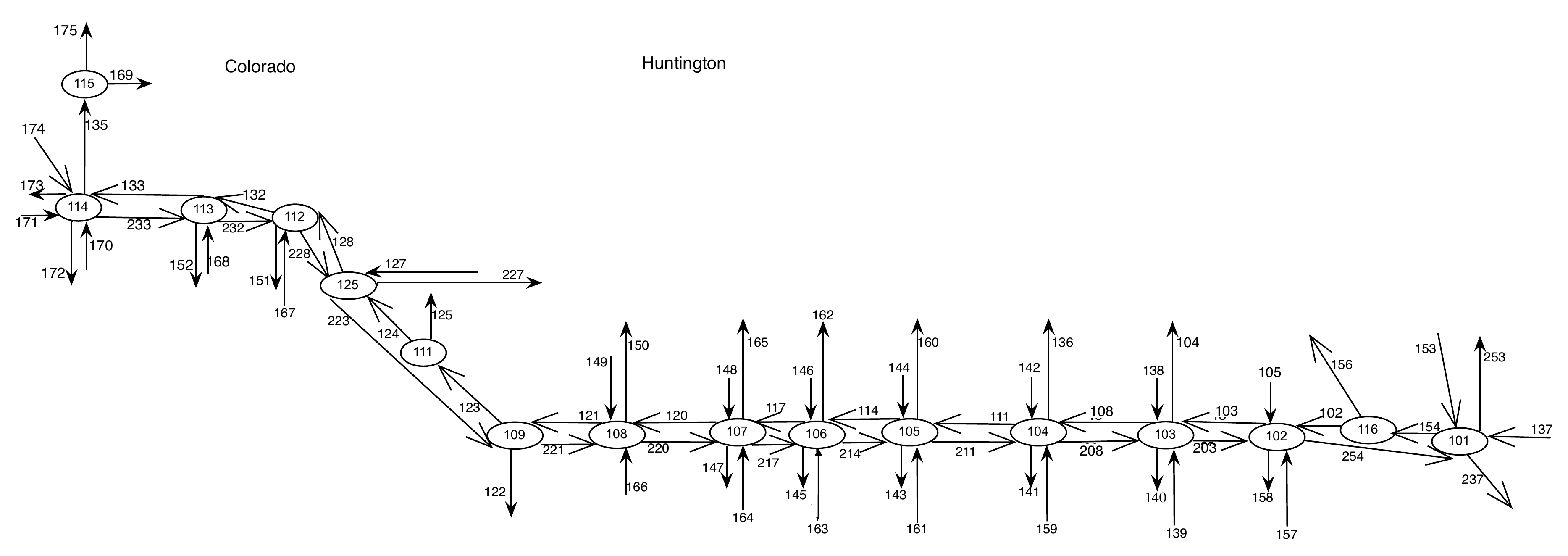}
\caption{Map of site and corresponding network graph. }
\label{fig-site}
\end{figure}

Figure \ref{fig-site} shows a map of the study site and its representation as a directed graph with 16 signalized intersections, 73 links, and 106 turn movements, hence 106 queues.  Each queue corresponds to a movement, so it is convenient to index a queue by a pair $(m,n)$ in which $m$ ($n$) is the incoming (outgoing)  link index.  For example, referring to Figure \ref{fig-site}, $x(139,104)$ is the number of vehicles in link 139 that are queued up at intersection 103 waiting to go to link 104.

For the first set of experiments we increase all the saturated flow rates  and the demands by the same factor $\G $ = 1.0, 1.5, 2.0, 2.5, and 3.0, with $\G=1.0$ being the base case.  For the second  set of experiments we replace the fixed-time control by the max pressure (MP) adaptive control under two switching regimes: MP4 permits four phase changes per cycle just like for  FT control, while MP6 permits six phase changes.  For each experiment  we present mean values of queue lengths and queuing delay.  Each stochastic simulation lasts 3 hours or 10,800 sec.  Since we want to calculate the mean values of queue lengths and delays, and since the queue length process is positive recurrent, it is reasonable to assume that the time average of these quantities over a three-hour long sample path will be  close to their statistical mean.  Hence the mean values reported below are the empirical time averages.
The results confirm all three predictions.

\textbf{Queue lengths} 
\begin{figure}[h!]
\centering
\includegraphics[width=5.8in]{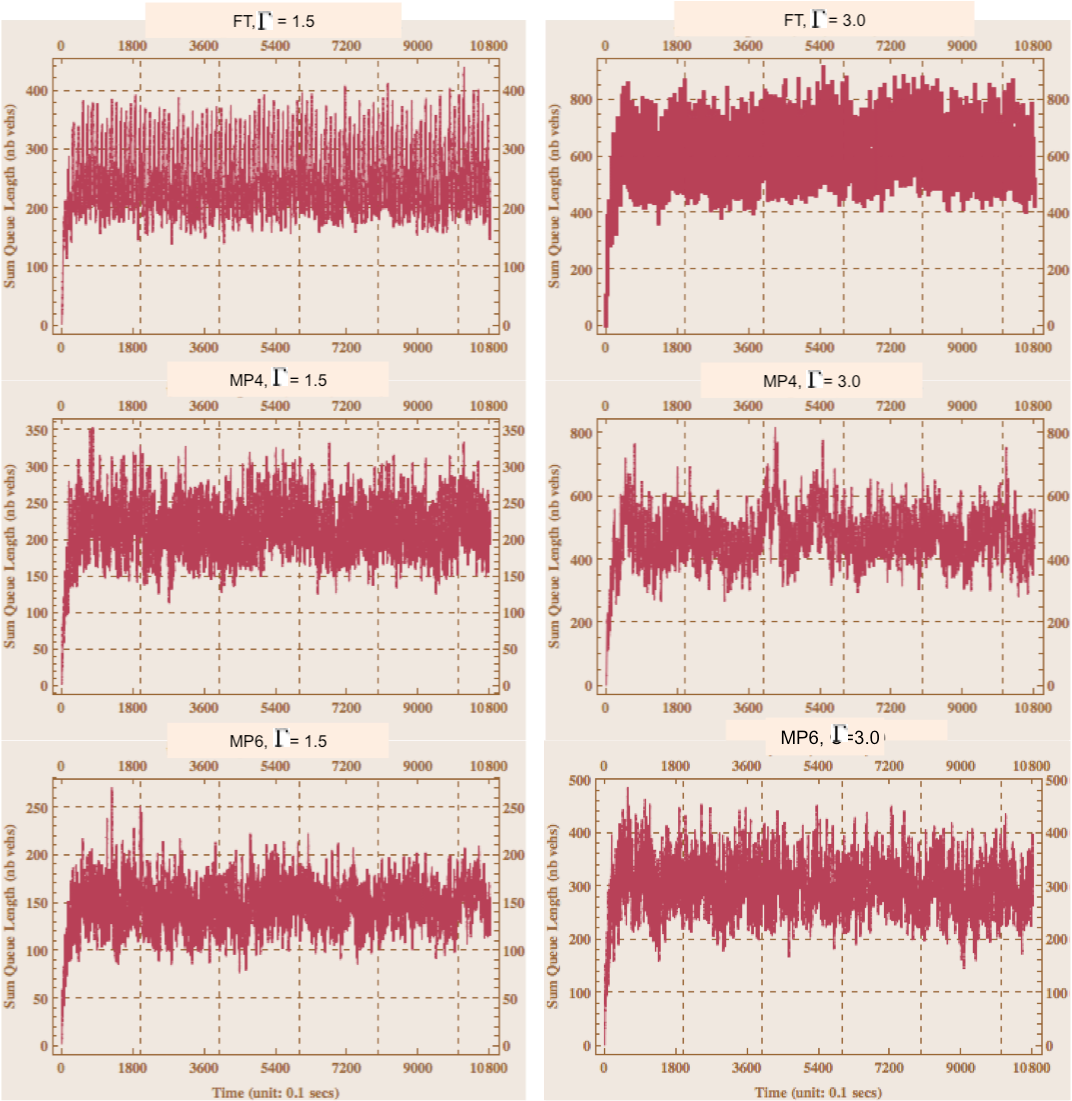}
\caption{Sum of all queues for FT, MP4, MP6 control when demand and saturation flows are scaled by $\G$ = 1.5 and 3.0.}
\label{sum_queues}
\end{figure}
Figure \ref{sum_queues} shows plots of the sum of all queues for three controls: FT (fixed time), MP4 (4 switches per cycle) and MP6 (6 switches per cycle) when the demand and saturation flow rates for the base case are scaled by $\G =1.5$ and 3.  Notice that the queue lengths for $\G=3$ are approximately twice as large as for $\G = 1.5$, for all three signal control schemes, as predicted by \eqref{14a}.

\begin{table}[h]
 \centering
 \begin{tabular}{|c|c|c|c|}
  \hline
 \bf $\G$, Total vph & \bf Signal Control  & \bf Mean Sum & \bf Ratio of sum to\\
 &\bf Type& \bf all queues (veh) &  \bf case $\G=1$\\ \hline 
 1.0, 14350&FT&197.4& 1\\
 1.0, 14350&MP4&146.5 & 1\\
 1.0, 14350&MP6&101.9 & 1\\
\hline
1.5, 21455&FT  & 248.6 & 1.26\\
1.5, 21455&MP4 & 213.7& 1.49 \\
1.5, 21455&MP6& 146.8 &1.46\\
 \hline
2.0, 28540&FT & 347.0 &1.76\\
2.0, 28540&MP4& 301.1  & 2.1\\
2.0, 28540&MP6 & 201.1 &1.99 \\
\hline
2.5, 35950&FT  & 455.7& 2.31 \\
2.5, 35950&MP4 & 395.0 & 2.7\\
2.5, 35950&MP6 &  248.7 & 2.46\\
\hline
3.0, 43298&FT& 586.2 &2.97\\
3.0, 43298&MP4 &  470.0 & 3.2\\
3.0, 43298&MP6 &  298.3& 2.95\\
\hline
\end{tabular}
\caption{Mean total queue length for different demands and signal control.}\label{table1}
\end{table}

Table \ref{table1} reports the average sum of all  queue lengths.  The first column gives the scale factor $\G$ and the average total exogenous demand in vph, e.g.\ the first column entry in the first row, (1.0, 14,350), refers to the base case $\G=1$ with an exogenous demand of 14,350 vph.  The second column indicates the signal control strategy used: FT is fixed time, MP4 is max pressure with 4 changes per cycle, and MP6 is max pressure with 6 changes per cycle.  The third column is the average number of vehicles summed over all 106 queues.  The fourth column is the ratio of the sum of queue lengths to the sum for the base case, $\G =1$.  As is seen in the fourth column, this ratio is roughly equal to the scale factor, $\G$.   For example for FT, the sum of queue lengths grows in the proportion 1: 1.26: 1.76: 2.31: 2.97 as the scale increases in the proportion 1: 1.5: 2.0: 2.5: 3.0.  This is in conformity with the first two predictions.  First, the queue sums are staying bounded, so the increase in the saturation flows by $\G$ does support demands that increase in the same proportion, and the queue lengths grow in the same proportion.   

\textbf{Queuing delay} We now consider queuing delay.  The prediction is that the mean queuing delay experienced by a vehicle stays the same despite the increase in demand.  Since there are 106 queues in all, we narrow our focus to the six queues at intersection 103.  Table \ref{table2} gives the delays of each queue as a function of five values of the scale factor $\G$ and three different signal control strategies: FT, MP4 and MP6.  Consider for example the queue $x(138,140)$.  As the saturation flow rates and the demand increase in proportion 1: 1.5: 2.0: 2.5: 3.0, the mean delay (sec) in this queue  under FT changes within a narrow range 22.35: 17.50: 18.99: 19.95: 20.69, while under MP4 the delay varies within the range [10.12, 11.30] and under MP6 the range is [6.73, 7.3].  When the delay is averaged over all queues at intersection 103, the variation with changes in demand is even smaller as is seen in Table \ref{table3}.
Note that three queues at intersection 103 correspond to right turn (RT) movements.  Since right turns on red are permitted the impact of increased rates of demand and saturation flow may not be adequately captured by the three models.  For this reason, Table \ref{table3} reports delays with and without including RT queues.

\begin{table}[h!]
 \centering
 \begin{tabular} {|c|c|c|r|r|r|}
 \hline
\bf $\G$ & \bf RT &\bf Queue & \bf Delay (sec) & \bf Delay (sec) & \bf Delay (sec) \\
 &  &\bf (Movement) & \bf FT & \bf MP4 & \bf MP6\\

\hline\hline
1.0&RT&$x$(103,104) &5.82&3.60 & 2.91\\
1.0&-&$x$(103,108) &3.07&13.60&8.94\\
1.0&RT&$x$(138,108) & 3.76&3.48&3.28\\
1.0&-&$x$(138,140) & 22.35& 11.30&7.30\\
1.0&RT&$x$(139,203) &2.48& 2.42&2.33\\
1.0&-&$x$(139,104) & 21.62& 11.17&7.49\\
 \hline\hline

1.5&RT&$x$(103,104) & 4.22&3.41 &2.46\\
1.5&-&$x$(103,108) &2.36 &13.70 & 7.14\\
1.5&RT&$x$(138,108) &0.80 &2.34 & 2.31\\
1.5&-&$x$(138,140) &17.50 & 10.71& 6.76\\
1.5&RT&$x$(139,203) & 0.73&1.57& 1.56\\
1.5&-&$x$(139,104) & 18.04 & 10.72& 6.77\\
 \hline\hline

2.0&RT&$x$(103,104) &5.08 &3.19 &2.27 \\
2.0&-&$x$(103,108) & 2.30& 14.28  &6.84  \\
2.0&RT&$x$(138,108) &0.88 &1.99&1.93\\
2.0&-&$x$(138,140) & 18.99&10.78&7.20\\
2.0&RT&$x$(139,203) & 0.75 &1.25&1.25\\
2.0&-&$x$(139,104) & 19.1&11.09&7.16\\
 \hline\hline

2.5&RT&$x$(103,104) &5.67 &3.48&2.28\\
2.5&-&$x$(103,108) &  2.19&14.91&6.77\\
2.5&RT&$x$(138,108)& 0.94&1.48&1.43\\
2.5&-&$x$(138,140) & 19.95&10.35&6.96\\
2.5&RT&$x$(139,203) & 0.76& 1.02&1.00\\
2.5&-&$x$(139,104) & 19.73&10.93&7.19\\
 \hline\hline

3.0&RT&$x$(103,104) &6.15 &  3.65  &2.30\\
3.0&-&$x$(103,108) &2.09 & 15.79  & 8.06 \\
3.0&RT&$x$(138,108) &1.08 & 1.08  &1.07 \\
3.0&-&$x$ (138,140) &20.69&  10.12  &6.73 \\
3.0&RT& $x$(139,203) &0.79& 0.79   &0.78 \\
3.0&-& $x$(139,104) &20.5&  10.13  &6.85 \\
 \hline\hline
\end{tabular}
\caption{Delays in all queues at intersection 103.  RT means right turn.}\label{table2}
\end{table}

\begin{table}[h!]
 \centering
 \begin{tabular}{|l||l||l|l|l|l|l||}
 \hline
 \bf $\G$ &\bf Delay & \bf Delay & \bf Delay& \bf Delay &\bf Delay &\bf Delay\\
& \bf w RT & \bf w RT & \bf w RT& \bf wo RT &\bf wo RT &\bf wo RT\\
& \bf FT &\bf MP4 & \bf MP6 &\bf FT & \bf MP4 & \bf MP6\\
\hline
1.0 & 9.9 &7.6 &5.4 & 15.7& 12.0 &7.9\\
1.5 & 7.3 &7.1 &4.5 & 12.7& 11.7 &6.9\\
2.0 & 7.9 &7.1 &4.4 & 13.5& 12.0 &7.1\\
2.5 & 8.2 &7.0 &4.3 & 14.0& 12.1 &7.0\\
3 &8.6 &6.9 &4.3 & 14.4& 12 &7.2\\
 \hline
\end{tabular}
\caption{Average delay at intersection 103. RT means right turn.}\label{table3}
\end{table}

\textbf{Faster phase switching}\
We now come to the third prediction: queue lengths will decline as the number of phase switches per cycle increases.   This is borne out when we compare the sum of queue lengths for MP4 vs MP6, for each level of demand.  MP4 permits four and MP6 permits 6 switches per cycle.  As Table \ref{table1} shows, the queues for MP6 are indeed smaller than for MP4. 

In fact the second and third models suggest a quantitative prediction.  In \eqref{11a} and \eqref{15a} $g$ is the switching speed up, $g = 1$, being the base case.  If we take MP4 as the base case, then MP6 corresponds to $g = 1.5$ (six vs four phase switches per cycle).  So according to   \eqref{11a} and \eqref{15a} the mean queue length under MP4 should be 1.5 times the mean queue length under MP6, for every demand.  Going back to Table \ref{table1}, we see that the ratios of the sum of queue lengths under MP4 to the length under MP are:
1.45, 1.46, 1.5, 1.6, and 1.6 for $\G = 1,1.5, 2.0, 2.5$ and 3.0. 
A similar prediction is upheld in Table \ref{table3}: the ratio of delays under MP4 and MP6 at queues in intersection 103 lie within a narrow band around 1.7 which can be compared to a switching speed up of $g = 1.5$.

\section{Conclusion} \label{sec-conc}
 
Intersections are the bottlenecks of urban roads, since their capacity is about one quarter of the maximum vehicle flow that can be accommodated by the approaches to a standard four-legged intersection.  This bottleneck capacity can be increased by a factor of two to three if vehicles are organized to cross the intersection in platoons with 0.75s headway at 45 mph or 0.7s headway at 30 mph to achieve a saturation flow rate of 4800 vph per movement.  
To organize themselves into a platoon when crossing an intersection, vehicles must have ACC capabilities that work when the vehicle is stopped, so that 
when the signal is actuated (light turns green) all queued vehicles  move together with a short headway.  All the drivers need to do is to engage the ACC
mode when stopped at an intersection. Under ACC control newly  arriving vehicles may join a moving platoon.

CACC capability can provide even shorter headway than ACC, because of its superior car-following performance.  CACC brings additional features, e.g. CACC vehicles may permit lane changing by a vehicle in a platoon.
Over the past five years, several platooning experiments on freeways have been reported in which factory-equipped ACC vehicles have been augmented by V2V communication capability.  The augmented vehicles can operate in CACC mode.  For example, \cite{ploeg} report a 6-vehicle CACC platoon, with a 0.5s headway; and \cite{cacc-milanes} describe  a 4-vehicle CACC platoon, with a 0.8s headway, capable of cut-in, cut-out and other maneuvers.
 There are additional questions that need  to be addressed when considering CACC in an urban road setting.   Should the intersection provide a target speed?
Should it suggest a minimum or maximum platoon size?  Should vehicles continue in platoon formation after crossing the intersection?

The paper explores the urban mobility benefits of larger saturation flow rates using queuing analysis and a simulation case study.   It reaches the following conclusions.
If the saturation flow rate is increased by a factor $\G$, the network can support an increase in demand by the same factor $\G$, with no increase in queuing delay or travel time, and using the same signal control.  However, the queues will also grow by the same factor $\G$, so if this leads to a saturation of the links, the improvement in throughput will be sub-linear in $\G$.  On the other hand, if the cycle time is reduced, the queues  will also be reduced, and this may restore the linear growth in demand.  Experimental results of platoons of CACC vehicles suggest that $\G = 2$ or 3 is technically achievable, which indicates a pure productivity increase of the urban road infrastructure by an unprecedented 200 to 300 percent using connected vehicle technology.  
The productivity increase can be shared between increased demand and a reduction in queuing delay by shortening the cycle time.  It is difficult to imagine any
technology in urban transportation that has such a great potential impact on mobility.

The analysis above has two significant limitations.  First, it only considers the limiting case in which there is 100 percent penetration of (C)ACC technology.  If 100 percent
penetration leads to a saturation flow gain of  $\G$, and a $p$ percent penetration rate leads to a smaller gain of  $\sigma (p) \times \G$, then
the analysis above will be valid with $\G$ replaced by $\sigma (p) \G$.  But it is a limitation of this study that we don't know how small $\sigma (p)$ is.
In a forthcoming paper we argue that for ACC, $\sigma (p) \approx p$, that is, the productivity gain is directly proportional to the penetration.  This  important result says that gains will occur immediately with the introduction of ACC.  However, for CACC, the gains are the same as for ACC for penetration rates $p$ below 50 percent, and begin to grow super-linearly in $p$  only for $p > 0.5 $.  The second limitation is that in short urban links vehicles
will slow down quickly as queues build up.  As a result the saturation flow rate at the upstream intersection will be reduced, thereby depriving the system of the full productivity benefit.  It is important to investigate  this reduction.

\bibliography{/Users/varaiya/Dropbox/varaiya-Main/varaiya/Bib/traffic}

\end{document}